\documentclass[12pt]{article}

\usepackage{mathrsfs}

\usepackage{amsthm,amsmath,amssymb,latexsym,color,cite}
\newtheorem{thm}{Theorem}[section]

\newtheorem{cor}[thm]{Corollary}
\newtheorem{conj}[thm]{Conjecture}
\newtheorem{prop}[thm]{Proposition}
\newtheorem{prob}[thm]{Problem}

\theoremstyle{remark}

\numberwithin{equation}{section}

\def\w{{\rm w}}

\def\i{{\rm i}}

\begin{document}

\begin{center}
{\Large\bf Pairs of lattice paths and positive trigonometric sums 
}
\end{center}

\vskip 2mm \centerline{Victor J. W. Guo$^1$  and Jiang Zeng$^{2}$}

\begin{center}
{\footnotesize $^1$Department of Mathematics, East China Normal
University,\\ Shanghai 200062,
 People's Republic of China\\
{\tt jwguo@math.ecnu.edu.cn,\quad http://math.ecnu.edu.cn/\textasciitilde{jwguo}}\\[10pt]
$^2$Universit\'e de Lyon; Universit\'e Lyon 1; Institut Camille
Jordan, UMR 5208 du CNRS;\\ 43, boulevard du 11 novembre 1918,
F-69622 Villeurbanne Cedex, France\\
{\tt zeng@math.univ-lyon1.fr,\quad
http://math.univ-lyon1.fr/\textasciitilde{zeng}} }
\end{center}

\vskip 0.7cm \noindent{\small{\bf Abstract.}
Ismail et al. (Constr. Approx. {\bf 15}  (1999) 69--81) proved the positivity of some trigonometric polynomials
with single binomial coefficients. In this paper, we prove
some similar results by replacing the  binomial coefficients with
products of two binomial coefficients. }

\vskip 0.5cm
\noindent{\bf MR Subject Classifications}: Primary: 42A32; Secondary: 05A15.
\section{Introduction}
Motivated by Bressoud's generalization of Borwein's conjectures (see \cite{An,Br,BW}),
Ismail, Kim and Stanton~\cite{IKS} proved that
the trigonometric polynomial
\begin{align}
\sum_{l}{M+N\choose M-kl} \cos(lx)
\label{eq:iks}
\end{align}
is a polynomial in $1+\cos(x)$ with nonnegative integral coefficients if $M,N$ and $k$ are positive integers satisfying
 $|M-N|\leq k$.
The starting point of \cite{IKS} is the observation that
the number of lattice paths from $O$ to $(M,N)$ not touching the lines $y=x\pm k$ is given by (see \cite[p. 12]{Na})
\begin{align}\label{eq:Nar}
\sum_{l}{M+N\choose M-kl}(-1)^l \geq 0\quad \textrm{if}\quad |M-N|\leq k.
\end{align}
By counting pairs of lattice paths (see Propositions 2.1 and 2.2)  we notice that
\begin{align}
&\sum_{l}{M+N\choose M-kl}^2(-1)^l \geq 0 \quad\text{if}\quad |M-N|\leq k, \label{eq:MN-2}\\
&\sum_{l}{M+N\choose M-kl}{M+N\choose N-kl}(-1)^l \geq 0. \label{eq:MNMN}
\end{align}
This encourages us to study \eqref{eq:iks} by replacing the single
binomial coefficient with a product of
binomial coefficients.
The following is our main result.
\begin{thm}\label{thm:MN2MN}
Let $M,N\geq 0$. Then the following two trigonometric polynomials
\begin{align}
&\sum_{l}{M+N\choose M-kl}^2 \cos(lx)
\quad\text{with \quad $|M-N|\leq k$},\label{eq1}\\
&\sum_{l}{M+N\choose M-kl}{M+N\choose N-kl} \cos(lx) \label{eq2}
\end{align}
are polynomials in $1+\cos(x)$ with nonnegative integral coefficients
\end{thm}

Ismail et al. \cite{IKS} generalized their formula \eqref{eq:iks}
to arbitrary polynomials. In the last section, we shall give such generalizations
of Theorems \ref{thm:MNap} and \ref{thm:MNap2}. It is interesting to notice that
 all the formulas given in \cite[Section~2]{IKS} are valid \emph{mutatis mutandis}.

\medskip

Finally we recall  the ${}_2F_1(1)$ or Chu-Vandermonde summation formula (see \cite[p. 67]{AAR}):
\begin{align}\label{eq:cv}
{}_2F_1(-n,a;c;1)=\frac{(c-a)_n}{(c)_n},
\end{align}
where $(a)_k=a(a+1)\ldots (a+k-1)$ and
${}_2F_1(a,b;c;x)=\sum_{n=0}^\infty \frac{(a)_n(b)_n}{(c)_n}\frac{x^n}{n!}$.

\section{Proof of Theorem 1.1}
A lattice path in the plan consists of unit steps of two types: in the north and east directions.
It is convenient to encode lattice paths by
 \emph{binary words} on the alphabet $\{0,1\}$.
For any word $u$ on $\{0,1\}$,
let $|u|_i$ denote the number of occurrences of $i\in \{0,1\}$ in the word $u$, and $|u|$ the length of $u$.
 Designate by $\mathscr{M}_{m,n}$ the set of words $u$ on $\{0,1\}$
 such that $|u|_1=m$ and $|u|_0=n$. Clearly the cardinality of $\mathscr{M}_{m,n}$ is ${m+n\choose n}$.

Throughout of this section we shall assume that $|M-N|\leq k$.
 For any integer $l$ let $\mathcal W_l=\mathscr M_{M-kl,N+kl}\times\mathscr M_{N+kl,M-kl}$ and
 $\mathcal W:=\cup_{l=-\infty}^\infty \mathcal W_l$.

\begin{prop}\label{prop:MN-2}
 The number of bi-words  $(u, v)\in \mathcal W_0$ such that $-k< \sum_{j\leq
s}(u_j-v_j)<k$  for all $s\leq M+N$ is given by \eqref{eq:MN-2}.
\end{prop}
\begin{proof}
Since
the cardinality of $\mathcal W_l$ is ${M+N\choose M-kl}^2$,
 weighting each term in $\mathcal W_l$ by $(-1)^l$,
the total  weight  of  $\mathcal W$  is given by \eqref{eq:MN-2}.

Let $ \mathcal W'$ (resp. $\mathcal W_0'$) be the set of
elements $(u,v)$ in $\mathcal W$ (resp. $\mathcal W_0$) satisfying the condition
$-k< \sum_{j\leq s}(u_j-v_j)<k$ for all $s\leq M+N$.
We first show that  $\mathcal W'=\mathcal W_0'$. Indeed,
if $(u,v)\in \mathcal W_l\cap \mathcal  W'$, then $|u|_1-|v|_1=M-N-2kl$.
As $|M-N|\leq k$, the hypothesis
$-k< |u|_1-|v|_1<k$  implies that $l=0$.
It remains to define a sign-reversing killing involution on $\mathcal W\setminus \mathcal W'$,
which reduces \eqref{eq:MN-2} to the
cardinality of  $\mathcal W'$.

For any $(u,v)\in\mathcal W\setminus \mathcal W'$, since $u_j-v_j\in \{0, \pm1\}$,
there exists some $s>0$ such that $\sum_{j\leq s}(u_j-v_j)=\pm k$.
Picking the smallest such $s$ and exchanging the first $s$ letters in $u$ and $v$, we obtain
 $u'=v_1\cdots v_s u_{s+1}\cdots u_{M+N}$ and
$v'=u_1\cdots u_s v_{s+1}\cdots v_{M+N}$. It is clear that if $(u,v)\in\mathcal W_j$ then
$(u',v')\in\mathcal W_{j-1}\cup \mathcal W_{j+1}$ and they have the same $s$.
Thus $(u,v)\mapsto (u',v')$ is a sign-reversing involution on $\mathcal W\setminus \mathcal W'$.
\end{proof}

Similarly, we can prove the following result.
\begin{prop}
The number of bi-words  $(u, v)$  on $\{0,1\}$, such that
$|u|_1=|v|_1=M$, $|u|_0=|v|_0=N$, and $-k< \sum_{j\leq s}(u_j-v_j)<k$ for all $s\leq M+N$, is given by \eqref{eq:MNMN}.
\end{prop}

For any $(u,v)\in \mathcal W_l$, we have $|u|_1-|v|_1=M-N-2kl$. So
we can uniquely factorize $u$ and $v$ into a product of subwords
$u^{(1)},\ldots,u^{(s+1)}$ and $v^{(1)},\ldots,v^{(s+1)}$,
respectively, with the largest $s$, such that $|u^{(j)}|=|v^{(j)}|$,
$|u^{(j)}|_1-|v^{(j)}|_1=\pm k$ for all $1\leq j\leq s$ and $-k\leq
|u^{(s+1)}|_1-|v^{(s+1)}|_1\leq k$. If
$|u^{(s+1)}|_1-|v^{(s+1)}|_1=\pm k$, then we say that $(u,v)$ has
$(s+1)$ {\it $k$-segments}; otherwise $(u,v)$ has $s$ {\it
$k$-segments}. We call $(u^{(1)},\ldots,u^{(s+1)};
v^{(1)},\ldots,v^{(s+1)})$ the {\it $k$-factorization} of $(u,v)$.
In addition, if there are $p$ values $r$ such that $\sum_{j=1}^r
(|u^{(j)}|_1-|v^{(j)}|_1)=\pm k$, we say that $(u,v)$ is of  {\it
class $p$} and denote by $\mathcal W_l(p)$ the set of elements in
$\mathcal W_l$ of class $p$.  Define  the set
\begin{align}\label{eq:condition}
B_p:=\left\{(u,v)\in\mathcal W_0(p)\colon
-2k<\sum_{j\leq s}(u_j-v_j)<2k \;\textrm{for}\; s=1,\ldots, M+N\right\}.
\end{align}

\begin{thm}\label{thm:MNap}
For $M,N\geq 0$, we have
\begin{align}
\sum_{l}{M+N\choose M-kl}^2 \cos(lx)=\sum_{p\geq 0} b_p(1+\cos x)^p,
\label{eq:MNap}
\end{align}
where $b_p$ is  the cardinality of $B_p$.
\end{thm}
\begin{proof} Setting $x=\pi$ in \eqref{eq:MNap} we see that
$b_0\geq 0$  by Proposition~\ref{prop:MN-2}. Assume that $p\geq 1$.
Writing
$\cos(lx)=T_l(\cos x)$, where
\begin{align*}
T_l(x)=(-1)^l\sum_{k=0}^l \frac{(-l)_k (l)_k}{k!(1/2)_k}(1+x)^k 2^{-k}
\end{align*}
is the Chebyshev polynomial of the first
kind~\cite[p. 101]{AAR},
we see that the coefficient  $b_p$
in \eqref{eq:MNap} is equal to
\begin{align*}
b_p
&=\sum_{l}{M+N\choose M-kl}^2 \frac{(-l)_p (l)_p}{p!(1/2)_p}\frac{(-1)^l}{2^p}\\
&=\sum_{|l|\geq p}{M+N\choose M-kl}^2 {|l|+p\choose |l|-p}\frac{2^p|l| (-1)^{l-p}}{|l|+p}.
\end{align*}
Therefore $b_p$ is the weight function of the set $\mathcal
W=\cup_{|l|\geq p}\mathcal W_l$ of bi-words  on $\{0,1\}$ with each
term in $\mathcal W_l$  weighted by
$$
\w(l,p)={|l|+p\choose |l|-p}\frac{2^p|l|}{|l|+p}(-1)^{|l|-p}.
$$
As the weight is $0$ for $|l|<p$, we may assume that $|l|\geq p>0$.
We will show that many of the terms (with weights) in $\mathcal W$ are summed to be $0$.
Moreover, all the other terms (with weights) will have $|l|=p$, and are thus positive.

By symmetry, we may assume that $N\geq M$.
Let $(u^{(1)},\ldots,u^{(s)}; v^{(1)},\ldots,v^{(s)})$ be
the $k$-factorization of $(u,v)\in \mathcal W$.
If all the quantities $|u^{(j)}|_1-|v^{(j)}|_1$ ($1\leq j\leq s$) are equal,
we say that $(u,v)$ is a {\it good guy}, otherwise it is a \emph{bad guy}.
\begin{itemize}
\item Assume that $N=M$. Suppose that $(u,v)\in\mathcal W_l$ is a good guy.  Then $s=2|l|$.
For any $0\leq r\leq |l|-p$,
if we choose $r$ pairs $(u^{(j)},v^{(j)})$ and
exchange $v^{(j)}$ and $u^{(j)}$,
then we obtain the  $k$-factorization of
a term $(u',v')$ in $\mathcal W_{l-r}$ if $l>0$, or in $\mathcal W_{l+r}$ if $l<0$.
The total weight of $(u',v')$ obtained in such a way is
\begin{align}
\sum_{r=0}^{|l|-p}{2|l|\choose r}\w(|l|-r,p). \label{eq:2.2a}
\end{align}
If $\ell:=|l|>p$, this sum  can be written as
\begin{align*}
\frac{2^p\ell}{\ell+p}{\ell+p\choose 2p}
&\bigg({}_2F_1(-2\ell, -\ell+p; -\ell-p+1;1)\\
&+\frac{2(\ell-p)}{\ell+p-1}{}_2F_1(-2\ell+1, -\ell+p+1; -\ell-p+2;1)\bigg),
\end{align*}
which is easily seen to be zero from the $_2F_1(1)$ evaluation \eqref{eq:cv}.

Note that any bad guy  in $\cup_{|l|\geq p}\mathcal W_l$ can be obtained
from a good guy $(u,v)$ by exchanging some factors $u^{(j)}$ and $v^{(j)}$ in its
$k$-factorization. It remains only to consider the good guys
$(u,v)$ in $\mathcal W_{-p} \cup \mathcal W_{p}$.

Since $\w(\pm p,p)=2^{p-1}$, each term in $\mathcal W_{-p} \cup \mathcal W_{p}$
must be counted $2^{p-1}$ times. We now give a surjection from the set $B_p$  to the set of good guys in
$\mathcal W_{-p} \cup \mathcal W_{p}$ so that each term
in the latter has $2^{p-1}$ preimages. Suppose that $(u,v)\in\mathcal W_0$ has class $p$ and
$(u^{(1)},\ldots,u^{(s)}; v^{(1)},\ldots,v^{(s)})$ is its $k$-factorization.
Then we must have $s=2p$.
If  $|u^{(s)}|_1-|v^{(s)}|_1=\varepsilon k$, where
$\varepsilon=\pm 1$,
then exchanging $u^{(j)}$ and
$v^{(j)}$ for any $j$ such that
$|u^{(j)}|_1-|v^{(j)}|_1=\varepsilon k$,
we obtain
a good guy  $(\overline{u},\overline{v})$ in $\mathcal W_{\varepsilon p}$.
Conversely, for any good guy
$(\overline{u},\overline{v})\in\mathcal W_{-p}\cup\mathcal W_{p}$,
let $(\overline{u}^{(1)},\ldots,\overline{u}^{(2p)};
\overline{v}^{(1)},\ldots,\overline{v}^{(2p)})$
be its $k$-factorization. Then there are exactly $2^{p-1}$
preimages of $(\overline{u},\overline{v})$ obtained
by exchanging $\overline{u}^{(j_i)}$ and $\overline{v}^{(j_i)}$,
where $j_i=2i-1$ or $2i$ for all $i=1,\ldots, p-1$, and $j_p=2p$.
This completes the proof of the $N=M$ case.

\item Assume that  $N>M$. By definition  $M+k\geq N$. Let $(u,v)\in\mathcal W_{l}$
be a good guy. Then $(u,v)$ has $2|l|$ or $(2|l|-1)$ $k$-segments if $l<0$, and
$2|l|$ or $(2|l|+1)$ $k$-segments if $l\geq 0$.
We can again choose any $r$ of the pairs $(u^{(j)},v^{(j)})$ in the $k$-factorization of
$(u,v)$ and exchange $v^{(j)}$ and $u^{(j)}$.
The cases of $2|l|$ $k$-segments sum to zero as before.
What remain are the cases of $(2|l|-1)$ $k$-segments for $l<0$
and of $(2|l|+1)$ $k$-segments for $l\geq 0$. The appropriate sum of
weights is given as follows:
\begin{align}
& B_{|l|,p}:=\sum_{r=0}^{|l|-p}{2|l|-1\choose r}\w(|l|-r,p), \label{eq:2.2b}\\
& C_{|l|,p}:=\sum_{r=0}^{|l|-p}{2|l|+1\choose r}\w(|l|-r,p). \label{eq:2.2c}
\end{align}
It follows from \eqref{eq:2.2a} that
$$
B_{|l|,p}+C_{|l|-1,p}=\sum_{r=0}^{|l|-p}{2|l|\choose r}\w(|l|-r,p)=0.
$$
Hence the total weight of
$\mathcal W_{-l-1}$ cancels that of $\mathcal W_{l}$.
One can check that for the maximal $l$ ($l=\lfloor M/k\rfloor$),
some terms in $\mathcal W_l$ have $s=2l+1$ if and only if
some terms in $\mathcal W_{-l-1}$ have $s=2l+1$. Thus this boundary
case is also canceled. Finally, the only remaining terms are
the:
\begin{itemize}
\item[(a)] good guys  in $\mathcal W_{-p}$;
\item[(b)] good guys in $\mathcal W_{p}$ with  $2p$ $k$-segments.
\end{itemize}
Again each term in $\mathcal W_{-p}\cup \mathcal W_p$ has
weight $2^{p-1}$.
We now give a surjection from the set $B_p$ to the set of good guys  in (a) and (b) by
distinguishing  the following two cases:
\begin{itemize}
\item If $N=M+k$, then (b) is empty. Suppose that $(u,v)\in B_p$
and\break $(u^{(1)},\ldots,u^{(2p-1)}; v^{(1)},\ldots,v^{(2p-1)})$ is its $k$-factorization.
Then $|u^{(2p-1)}|-|v^{(2p-1)}|\break =-k$, and for any $i=1,\ldots, p-1$,
exactly one of $|u^{(2i-1)}|-|v^{(2i-1)}|$ and $|u^{(2i)}|-|v^{(2i)}|$ is
equal to $-k$. For all $j$, exchanging $u^{(j)}$ and $v^{(j)}$ if $|u^{(j)}|_1-|v^{(j)}|_1=-k$,
we obtain a good guy  in $\mathcal W_{-p}$.
\item If $M+k>N>M$, for any $(u,v)\in B_p$, let
$(u^{(1)},\ldots,u^{(s)}; v^{(1)},\ldots,v^{(s)})$
 be
its $k$-factorization. We have $s=2p+1$ or $s=2p$.
If $|u^{(s-1)}|_1-|v^{(s-1)}|_1=-k$, then for all $j$,
exchanging $u^{(j)}$ and $v^{(j)}$ if $|u^{(j)}|_1-|v^{(j)}|_1=-k$,
we obtain a good guy in $\mathcal W_{-p}$.
If $|u^{(s-1)}|_1-|v^{(s-1)}|_1=k$ (in this case we must have $s=2p+1$), then
exchanging $u^{(j)}$ and $v^{(j)}$ if $|u^{(j)}|_1-|v^{(j)}|_1=k$,
we obtain a good guy  in $\mathcal W_{p}$ with  $2p$ $k$-segments.
\end{itemize}
It is easy to see that each term in (a) and (b) has $2^{p-1}$ preimages.
\end{itemize}
This completes the whole proof.
\end{proof}

Theorem~\ref{thm:MNap} has the following  sister theorem.
\begin{thm}\label{thm:MNap2}
For $M,N\geq 0$, we have
$$
\sum_{l}{M+N\choose M-kl}{M+N\choose N-kl} \cos(lx)=\sum_{p\geq 0} c_p(1+\cos(x))^p,
$$
where $c_p$ is the number of $(u,v)\in\mathscr{M}_{M,N}\times \mathscr{M}_{M,N}$ of class $p$ such that
$-2k< \sum_{j\leq s}(u_j-v_j)<2k$ for all $s=1,\ldots,M+N$.
\end{thm}

For $k=1$ and $N=M$, since the number of good guys $(u, v)$ on
$\{0,1\}$ such that $|u|=|v|=2M$ and $|u|_1=|v|_0=M+p$ is given by
${2M\choose M+p}{M+p\choose M-p}$, we get the following explicit formula:
\begin{align*}
\sum_{l=-M}^M{2M\choose M-l}^2 \cos(lx)
=\sum_{p=0}^M\frac{(2M)!2^p}{(M-p)!^2(2p)!}(1+\cos(x))^p.
\end{align*}
Similarly, for $k=1$ and $N=M+1$, we have
\begin{align*}
\sum_{l=-M-1}^M{2M+1\choose M-l}^2 \cos(lx)
& =\sum_{p=1}^{M+1}\frac{(2M+1)!2^{p-1}}{(M-p+1)!^2(2p-1)!}(1+\cos(x))^p, \\
\sum_{l=-M}^M {2M+1\choose M-l}{2M+1\choose M+l} \cos(lx)
& =\sum_{p=0}^{M}\frac{(2M+1)!2^p}{(M-p)!(M-p+1)!(2p)!}(1+\cos(x))^p.
\end{align*}
\section{Further extensions}
In the proof of Theorem~\ref{thm:MNap}, the nonnegativity of the coefficient
of $(1+\cos(x))^p$ for $p>0$ follows from the nonnegativity of
\eqref{eq:2.2a}. Thus the proof of Theorem~\ref{thm:MNap} is also valid
for certain weights other than the $T$-Chebyshev weight.
\begin{thm}\label{thm:extension}
Suppose that $p_l(z)=\sum_{p=0}^l \w(l,p)z^p$ is a polynomial in $z$ of degree at most
$l$. If \eqref{eq:2.2a} is nonnegative, then
\begin{align*}
\sum_{l}{M+N\choose M-kl}^2 p_{|l|}(z)
&=\sum_{p\geq 0}b_p z^p \quad\text{{\rm(}for $|M-N|\leq k${\rm)}},\\
\sum_{l}{M+N\choose M-kl}{M+N\choose N-kl} p_{|l|}(z)
&=\sum_{p\geq 0}c_p z^p,
\end{align*}
where $b_p\geq 0$ and $c_p\geq 0$ for $p> 0$.
\end{thm}

For example, if
$$
p_l(z)=\frac{(\alpha+\beta+1)_l}{(\beta+1)_l}P_{l}^{(\alpha,\beta)}(z-1)
$$
where $P_n^{(\alpha,\beta)}(x)$ is
the Jacobi polynomials defined by
$$
P_n^{(\alpha,\beta)}(x)=\frac{(\alpha+1)_n}{n!}{}_2F_1(-n,n+\alpha+\beta+1; \alpha+1;(1-x)/2)
$$
(see \cite[p. 99]{AAR}), then the $_2F_1(1)$ evaluation  \eqref{eq:cv} yields that
\begin{align}
\eqref{eq:2.2a}=\frac{(l+\alpha+\beta+1)_p (\alpha+\beta+1)_l (l-p-\alpha-\beta)_{l-p}}
{(l-p)!p!(\beta+1)_p (l+\alpha+\beta+2p)_{l-p}2^p}.
\label{eq:jacobi}
\end{align}

It is clear that if $-1\leq \alpha+\beta\leq 1$ and $\beta>-1$, we have
$\eqref{eq:jacobi}\geq 0$.
\begin{cor}Let $-1\leq \alpha+\beta\leq 1$ and $\beta>-1$. Then
\begin{align*}
\sum_{l}{M+N\choose M-kl}^2
\frac{(\alpha+\beta+1)_{|l|}}{(\beta+1)_{|l|}} P_{|l|}^{(\alpha,\beta)}(z-1)
&=\sum_{p\geq 0}b_p z^p \quad\text{if $|M-N|\leq k$}, \\
\sum_{l}{M+N\choose M-kl}{M+N\choose N-kl}
\frac{(\alpha+\beta+1)_{|l|}}{(\beta+1)_{|l|}} P_{|l|}^{(\alpha,\beta)}(z-1)
&=\sum_{p\geq 0}c_p z^p,
\end{align*}
where $b_p\geq 0$ and $c_p\geq 0$ for $p>0$.
\end{cor}

There is a \emph{sine} version of Theorem~\ref{thm:MN2MN}.
By derivation,  Theorems~\ref{thm:MNap}
and \ref{thm:MNap2} imply that
\begin{align*}
\sum_{l}{M+N\choose M-kl}^2 \frac{l\sin(lx)}{\sin(x)}
&\geq 0 \quad\text{if\quad $|M-N|\leq k$},\\
\sum_{l}{M+N\choose M-kl}{M+N\choose N-kl} \frac{l\sin(lx)}{\sin(x)}
&\geq 0,
\end{align*}
for any real $x$. Another version will be given in Corollary~\ref{cor:mnsin}.

We now apply Theorem~\ref{thm:extension} to the polynomial
\begin{align*}
p_l(z)=
\begin{cases}
0,&\text{if $l=0$}, \\[5pt]
\displaystyle\frac{(\alpha+\beta+1)_{l-1}}{(\beta+1)_{l-1}}P_{l-1}^{(\alpha,\beta)}(z-1),
&\text{if $l>0$}.
\end{cases}
\end{align*}
The argument of Theorem \ref{thm:MNap} also implies that the constant term
is nonnegative in this case. To verify that Theorem~\ref{thm:extension} may
be used, we again apply the $_2F_1(1)$ evaluation to find that
\begin{align*}
\eqref{eq:2.2a}=\frac{(l+\alpha+\beta)_p (\alpha+\beta+1)_{l-1} (l-p-\alpha-\beta+1)_{l-p-1}}
{(l-p-1)!p!(\beta+1)_p (l+\alpha+\beta+2p-1)_{l-p-1}2^p},
\end{align*}
so that the nonnegativity holds if $-1\leq\alpha+\beta\leq 2$ and $\beta>-1$. For
the special case $\alpha=\beta=1/2$, we obtain the Chebyshev polynomials of the second kind,
thus the next corollary.
\begin{cor}\label{cor:mnsin}
Let $M,N\geq 0$. Then
\begin{align*}
\sum_{l}{M+N\choose M-kl}^2 \frac{\sin(|l|x)}{\sin(x)}
&\geq 0 \quad\text{if\quad $|M-N|\leq k$},\\
\sum_{l}{M+N\choose M-kl}{M+N\choose N-kl} \frac{\sin(|l|x)}{\sin(x)}
&\geq 0,
\end{align*}
for any real $x$.
\end{cor}

It seems that Theorem \ref{thm:MN2MN} can be further generalized. For example, we make the following conjecture.
\begin{conj}\label{conj:mn}
Let $M_i$ and $ N_i$ be nonnegative integers such that
$|M_i-N_i|\leq k$ for $i=1,\ldots, r$. Then
\begin{align} \label{conj1}
\sum_{l}\prod_{i=1}^r {M_i+N_i\choose M_i-kl}\cos(lx)
\end{align}
is a polynomial in $1+\cos(x)$ with nonnegative integral coefficients, and
\begin{align*}
\sum_{l}\prod_{i=1}^r {M_i+N_i\choose M_i-kl}\frac{\sin(|l|x)}{\sin(x)}\geq 0
\end{align*}
for any real $x$.
\end{conj}

\medskip
\noindent\emph{Remark.}
Jaming Philippe (personal communication)
has observed  that the nonnegativity  of \eqref{conj1} for any real $x$
can be derived from the $r=1$  case as follows.
For any two trigonometric polynomials
$P(x)=\sum_k a_ke^{\i kx}$ and $Q(x)=\sum_k b_ke^{\i kx}$,
we have the usual convolution formula:
$$
(P*Q)(x):=\sum_k a_kb_ke^{\i kx}=\frac{1}{2\pi}\int_{-\pi}^\pi P(t)Q(x-t)dt,
$$
where $\i=\sqrt{-1}$. For the cosine polynomials, we just have the further
constraint $a_{-k}=a_k$. Clearly, the convolution theorem implies that if $P$ and $Q$ are nonnegative polynomials, then
so is $P*Q$.  Moreover, it
obviously extends to the product of $r\geq 3$ polynomials:
\begin{align} \label{eq:new}
\sum_l \prod_{i=1}^r {M_i+N_i\choose M_i-k_i}
\cos(lx)
=(P_1*P_2*\cdots*P_N)(x),
\end{align}
where
$$
P_i(x)=\sum_l {M_i+N_i\choose M_i-k_i}\cos(lx).
$$
For $|M_i-N_i|\leq k_i$, $i=1,\ldots,r$, the nonnegativity of \eqref{eq:new} follows from that of
\eqref{eq:iks}. Nevertheless, our Conjecture \ref{conj:mn} still remains challenging.

When $M_i=N_i$ ($i=1,\ldots, r$), $x=\pi$ and $k=1$ the nonnegativity of \eqref{conj1}
 was also proved in \cite{GHZ} along with a $q$-analogue.
   Some $q$-analogues of
 \eqref{eq:Nar}  were considered in \cite{ABBBFV,Br,IKS}. We  would like to end this paper with the following
 problem.
\begin{prob}
Find a $q$-analogue of Theorem {\rm\ref{thm:MN2MN}}.
\end{prob}

\medskip
\noindent{\bf Acknowledgments.} This work was done during the first author's visit to
Institut Camille Jordan of Univerit\'e Lyon I, and was supported by
Project MIRA 2007 de la R\'egion Rh\^one-Alpes.
The first author was also sponsored by Shanghai Educational Development Foundation under the Chenguang Project 2007CG29
and Shanghai Leading Academic Discipline Project, Project Number: B407.
\renewcommand{\baselinestretch}{1}

\end{document}